\documentclass{article}

\usepackage[english]{babel}

\usepackage[letterpaper,top=2cm,bottom=2cm,left=3cm,right=3cm,marginparwidth=1.75cm]{geometry}

\usepackage{amsmath}
\usepackage{amssymb}
\usepackage{amsthm}
\usepackage{tabstackengine}
\stackMath
\makeatletter
\renewcommand\TAB@delim[1]{\scriptstyle#1}
\makeatother
\setstackgap{S}{2pt}
\usepackage{xcolor}
\usepackage{graphicx}
\usepackage[colorlinks=true, hidelinks]{hyperref}
\newtheorem{theorem}{Theorem}
\newtheorem{lemma}{Lemma}

\begin{document}
\begin{center}
		\vskip 1cm{\LARGE \bf Closed-Form Evaluation of Two Apéry-Like Series of Weight $4$}\\		
		\vskip 5mm
		Jorge Antonio González Layja\footnotemark\\
        Mexico
	\end{center}
\footnotetext{Email: \href{mailto:jorgelayja16@gmail.com}{jorgelayja16@gmail.com}}
\begin{abstract}
This paper presents closed-form evaluations of two new Apéry-like series of weight $4$ that involve harmonic numbers of the form $H_{2k}$. Several key results are derived and subsequently used to establish connections to the main series.
\end{abstract}

\\[1ex]

\section{Introduction and Preliminaries}
The principal aim of this paper is to obtain closed-form evaluations of the following series:
$$\sum _{k=1}^{\infty }\frac{4^kH_{2k}^{\left(2\right)}}{k^2\binom{2k}{k}}\quad \operatorname{and}\quad\sum _{k=1}^{\infty }\frac{4^kH_{2k}^2}{k^2\binom{2k}{k}},$$
where $\binom{2k}{k}$ denotes the central binomial coefficient. To the best of our knowledge, these series have not previously appeared in the literature in closed form.
\\[1.5ex]
The $k$-th generalized harmonic number of order $m$, denoted by $H_k^{\left(m\right)}$, is defined as
$$H_k^{\left(m\right)}=\sum _{n=1}^k\frac{1}{n^m},\quad m\in \mathbb{Z}^+,$$
where $H_k^{\left(1\right)}=H_k$ is the $k$-th (ordinary) harmonic number.
\\[1.5ex]
An Apéry-like series involving harmonic numbers, of weight $m_1+m_2+\cdots +m_n+n$, is given by the summation
$$\sum _{k=1}^{\infty }\frac{4^kH_k^{\left(m_1\right)}H_k^{\left(m_2\right)}\cdots H_k^{\left(m_n\right)}}{k^n\binom{2k}{k}},$$
where $m_1,m_2,\dots ,m_n$ and $n$ are positive integers. According to this definition, we deduce that each of the two main series has weight $4$.
\\[1.5ex]
The Riemann zeta function is defined as
$$\zeta \left(s\right)=\sum _{k=1}^{\infty }\frac{1}{k^s},\quad \mathfrak{R} \left(s\right)>1.$$
The polylogarithm of order $s$ is defined as
$$\operatorname{Li}_s\left(z\right)=\sum _{k=1}^{\infty }\frac{z^k}{k^s},\quad \left|z\right|< 1.$$
Furthermore, if $\mathfrak{R} \left(s\right)>1$, the series converges on the closed disc $\left|z\right|\le 1$.
\\[1.5ex]
The inverse tangent integral of order $n\ge 2$ is defined as
\begin{equation*}
\operatorname{Ti}_n\left(z\right)=\sum _{k=1}^{\infty }\frac{\left(-1\right)^{k+1}z^{2k-1}}{\left(2k-1\right)^n},\quad \left|z\right|\le 1.\tag{1.1}
\label{1.1}
\end{equation*}
For the special case $n=1$, we have
\begin{equation*}
\operatorname{Ti}_1\left(z\right)=\sum _{k=1}^{\infty }\frac{\left(-1\right)^{k+1}z^{2k-1}}{2k-1}=\arctan \left(z\right),\quad \left|z\right|\le 1,\, z\neq \pm i.\tag{1.2}
\label{1.2}
\end{equation*}
If $n=2$ and $z=1$, we have
\begin{equation*}
\operatorname{Ti}_2\left(1\right)=\sum _{k=1}^{\infty }\frac{\left(-1\right)^{k+1}}{\left(2k-1\right)^2}=G,\tag{1.3}
\label{1.3}
\end{equation*}
where $G$ denotes Catalan's constant.
\\[1.5ex]
Additionally, the inverse tangent admits the following complex logarithmic representation:
$$\arctan \left(z\right)=-\frac{i}{2}\ln \left(\frac{\left(1+iz\right)^2}{1+z^2}\right),$$
for complex $z$ except along two vertical branch cuts extending upward and downward from $i$ and $-i$, respectively.
\\[1.5ex]
Moreover, this identity can be rearranged as
\begin{equation*}
\ln \left(\frac{1-z}{1-iz}\right)=\frac{1}{2}\ln \left(\frac{\left(1-z\right)^2}{1+z^2}\right)+i\arctan \left(z\right),\quad \operatorname{\mathfrak{R}} \left(z\right)<0.\tag{1.4}
\label{1.4}
\end{equation*}
\section{The Lemmas and their Proofs}

\begin{lemma} \label{lma1}The following equalities hold:
\begin{alignat*}{2}
\left(\operatorname{i}\right)&\quad&&\int _0^1x^{k-1}\ln ^2\left(1-x\right)\:dx=\frac{H_k^2+H_k^{\left(2\right)}}{k},\quad k\in \mathbb{Z}^+;\\
\left(\operatorname{ii}\right)&\quad&&\frac{1}{2}\sum _{k=1}^{\infty }\frac{4^kx^{2k-1}}{k\binom{2k}{k}}=\frac{\arcsin \left(x\right)}{\sqrt{1-x^2}},\quad \left|x\right|<1;\\
\left(\operatorname{iii}\right)&\quad&&\sum _{k=1}^{\infty }\frac{4^k}{k^2\binom{2k}{k}}=3\zeta \left(2\right);\\
\left(\operatorname{iv}\right)&\quad&&\sum _{k=1}^{\infty }\frac{4^kH_k^{\left(2\right)}}{k^2\binom{2k}{k}}=-\zeta \left(4\right)+8\operatorname{Li}_4\left(\frac{1}{2}\right)+4\ln ^2\left(2\right)\zeta \left(2\right)+\frac{1}{3}\ln ^4\left(2\right).
\end{alignat*}
\begin{proof}
The identity in point $\left(\operatorname{i}\right)$ is proved in \cite[p. 60]{book2}. Two proofs for the generating function in point $\left(\operatorname{ii}\right)$ are presented in \cite[pp. 331-333]{book2}. The series in point $\left(\operatorname{iii}\right)$ is obtained by integrating the previous generating function and evaluating at $x=1$. The result in point $\left(\operatorname{iv}\right)$ is demonstrated in \cite[pp. 367-369]{book1} and \cite[p. 334]{book2}.
\end{proof}
\end{lemma}

\begin{lemma} \label{lma2}The following equalities hold:
\begin{alignat*}{2}
\left(\operatorname{i}\right)&\quad&&\tan \left(x\right)\ln \left(\sin \left(x\right)\right)=-\sum _{k=1}^{\infty }\left(\int _0^1\frac{1-t}{1+t}t^{k-1}\:dt\right)\sin \left(2kx\right),\quad 0<x<\frac{\pi }{2};\\
\left(\operatorname{ii}\right)&\quad&&\int _0^{\frac{\pi }{2}}x^2\csc \left(x\right)\sin \left(2kx\right)\:dx=3\zeta \left(2\right)\sum _{n=1}^k\frac{\left(-1\right)^{n+1}}{2n-1}-4\sum _{n=1}^k\frac{\left(-1\right)^{n+1}}{\left(2n-1\right)^3},\quad k\in \mathbb{Z}_{\ge 0};\\
\left(\operatorname{iii}\right)&\quad&&\int _0^{\frac{\pi }{2}}x^2\sec \left(x\right)\ln \left(\sin \left(x\right)\right)\:dx=\frac{45}{16}\zeta \left(4\right)-4G^2.
\end{alignat*}
\begin{proof}
The identity in point $\left(\operatorname{i}\right)$ is proved in \cite[p. 243]{book2}.
\\[1.5ex]
To prove the identity in point $\left(\operatorname{ii}\right)$, we will first denote the integral in such point as $I_k$ and consider calculating $I_{n+1}-I_n$. This means that
$$I_{n+1}-I_n=\int _0^{\frac{\pi }{2}}x^2\csc \left(x\right)\left(\sin \left(2\left(n+1\right)x\right)-\sin \left(2nx\right)\right)\:dx=2\int _0^{\frac{\pi }{2}}x^2\cos \left(\left(2n+1\right)x\right)\:dx.$$
Hence, by applying integration by parts and noting that $\sin \left(\pi n\right)=0$ and $\cos \left(\pi n\right)=\left(-1\right)^n$ for $n\in \mathbb{Z}_{\ge 0}$, we obtain
$$I_{n+1}-I_n=3\zeta \left(2\right)\frac{\left(-1\right)^n}{2n+1}-4\frac{\left(-1\right)^n}{\left(2n+1\right)^3}.$$
Summing both sides from $n=0$ to $k-1$ yields a telescoping left-hand side. Thus,
$$I_k-I_0=3\zeta \left(2\right)\sum _{n=0}^{k-1}\frac{\left(-1\right)^n}{2n+1}-4\sum _{n=0}^{k-1}\frac{\left(-1\right)^n}{\left(2n+1\right)^3},$$
and, after reindexing the sums on the right-hand side and observing that $I_0=0$, we arrive at the desired result.
\\[1.5ex]
In order to calculate the integral in point $\left(\operatorname{iii}\right)$, let us consider the identity in point $\left(\operatorname{ii}\right)$. By multiplying its both sides by $-\frac{1-t}{1+t}t^{k-1}$, integrating from $t=0$ to $t=1$, and taking the sum from $k=1$ to $\infty$, we obtain
$$\int _0^{\frac{\pi }{2}}x^2\csc \left(x\right)\left(-\sum _{k=1}^{\infty }\left(\int _0^1\frac{1-t}{1+t}t^{k-1}\:dt\right)\sin \left(2kx\right)\right)\:dx$$
$$=-3\zeta \left(2\right)\int _0^1\frac{1-t}{1+t}\left(\sum _{k=1}^{\infty }t^{k-1}\sum _{n=1}^k\frac{\left(-1\right)^{n+1}}{2n-1}\right)\:dt+4\int _0^1\frac{1-t}{1+t}\left(\sum _{k=1}^{\infty }t^{k-1}\sum _{n=1}^k\frac{\left(-1\right)^{n+1}}{\left(2n-1\right)^3}\right)\:dt.$$
Furthermore, by applying the identity in point $\left(\operatorname{i}\right)$ to the left-hand side and reversing the order of summations on the right-hand side, we get
$$\int _0^{\frac{\pi }{2}}x^2\csc \left(x\right)\tan \left(x\right)\ln \left(\sin \left(x\right)\right)\:dx$$
$$=-3\zeta \left(2\right)\int _0^1\frac{1-t}{1+t}\left(\sum _{n=1}^{\infty }\frac{\left(-1\right)^{n+1}}{2n-1}\sum _{k=n}^{\infty }t^{k-1}\right)\:dt+4\int _0^1\frac{1-t}{1+t}\left(\sum _{n=1}^{\infty }\frac{\left(-1\right)^{n+1}}{\left(2n-1\right)^3}\sum _{k=n}^{\infty }t^{k-1}\right)\:dt,$$
and, since $\sum _{k=n}^{\infty }t^{k-1}=\frac{t^{n-1}}{1-t}$, we have
$$\int _0^{\frac{\pi }{2}}x^2\sec \left(x\right)\ln \left(\sin \left(x\right)\right)\:dx$$
$$=-3\zeta \left(2\right)\int _0^1\frac{1}{1+t}\left(\sum _{n=1}^{\infty }\frac{\left(-1\right)^{n+1}t^{n-1}}{2n-1}\right)\:dt+4\int _0^1\frac{1}{1+t}\left(\sum _{n=1}^{\infty }\frac{\left(-1\right)^{n+1}t^{n-1}}{\left(2n-1\right)^3}\right)\:dt.$$
Moreover, by applying the change of variables $t\mapsto t^2$ to both expressions and using the definitions in $\left(\ref{1.1}\right)$ and $\left(\ref{1.2}\right)$, we obtain
$$\int _0^{\frac{\pi }{2}}x^2\sec \left(x\right)\ln \left(\sin \left(x\right)\right)\:dx=-6\zeta \left(2\right)\int _0^1\frac{\arctan \left(t\right)}{1+t^2}\:dt+8\int _0^1\frac{\operatorname{Ti}_3\left(t\right)}{1+t^2}\:dt.$$
By applying integration by parts to the rightmost integral and noting that $\frac{d}{dt}\operatorname{Ti}_n\left(t\right)=\frac{\operatorname{Ti}_{n-1}\left(t\right)}{t}$, it follows that
$$\int _0^{\frac{\pi }{2}}x^2\sec \left(x\right)\ln \left(\sin \left(x\right)\right)\:dx$$
$$=-3\zeta \left(2\right) \arctan ^2\left(t\right)\Big|_{t=0}^{t=1}+8 \arctan \left(t\right)\operatorname{Ti}_3\left(t\right)\Big|_{t=0}^{t=1}-8\int _0^1\frac{\arctan \left(t\right)\operatorname{Ti}_2\left(t\right)}{t}\:dt$$
$$=-3\zeta \left(2\right)\arctan ^2\left(1\right)+8\arctan \left(1\right)\operatorname{Ti}_3\left(1\right)-4\operatorname{Ti}_2^2\left(1\right).$$
Therefore, by applying the definition in $\left(\ref{1.3}\right)$ together with the known values $\arctan \left(1\right)=\frac{\pi }{4}$ and $\operatorname{Ti}_3\left(1\right)=\frac{\pi ^3}{32}$, we conclude the proof.
\end{proof}
\end{lemma}

\begin{lemma} \label{lma3}The following equalities hold:
\begin{alignat*}{2}
\left(\operatorname{i}\right)&\quad&&\operatorname{\mathfrak{R}} \left\{\operatorname{Li}_4\left(1\pm i\right)\right\}=\frac{485}{512}\zeta \left(4\right)-\frac{5}{16}\operatorname{Li}_4\left(\frac{1}{2}\right)+\frac{1}{8}\ln ^2\left(2\right)\zeta \left(2\right)-\frac{5}{384}\ln ^4\left(2\right);\\
\left(\operatorname{ii}\right)&\quad&&\int _0^{\frac{\pi }{2}}x\ln ^2\left(1-\sin \left(x\right)\right)\:dx=\frac{125}{8}\zeta \left(4\right)-5\operatorname{Li}_4\left(\frac{1}{2}\right)+2\ln ^2\left(2\right)\zeta \left(2\right)-\frac{5}{24}\ln ^4\left(2\right).
\end{alignat*}
\begin{proof}
The identity in point $\left(\operatorname{i}\right)$ is proved in \cite[pp. 777-778]{book3}.
\\[1.5ex]
To evaluate the integral in point $\left(\operatorname{ii}\right)$, we begin by applying the change of variables $x\mapsto 2\arctan \left(x\right)$. This yields
\begin{equation*}
\int _0^{\frac{\pi }{2}}x\ln ^2\left(1-\sin \left(x\right)\right)\:dx=4\int _0^1\frac{\arctan \left(x\right)\ln ^2\left(\frac{\left(1-x\right)^2}{1+x^2}\right)}{1+x^2}\:dx.\tag{2.1}
\label{2.1}
\end{equation*}
Next, by raising both sides of $\left(\ref{1.4}\right)$ to the power of 3, replacing $z$ with $x$, and taking the imaginary part, we obtain
$$\operatorname{\mathfrak{I}} \left\{\ln ^3\left(\frac{1-x}{1-ix}\right)\right\}=\frac{3}{4}\arctan \left(x\right)\ln ^2\left(\frac{\left(1-x\right)^2}{1+x^2}\right)-\arctan ^3\left(x\right),\quad x\in \mathbb{R}_{<1},$$
and, if we isolate $\arctan \left(x\right)\ln ^2\left(\frac{\left(1-x\right)^2}{1+x^2}\right)$, substitute it into $\left(\ref{2.1}\right)$, and expand, we get
\begin{equation*}
\int _0^{\frac{\pi }{2}}x\ln ^2\left(1-\sin \left(x\right)\right)\:dx=\frac{16}{3}\int _0^1\frac{\arctan ^3\left(x\right)}{1+x^2}\:dx+\frac{16}{3}\operatorname{\mathfrak{I}} \left\{\int _0^1\frac{\ln ^3\left(\frac{1-x}{1-ix}\right)}{1+x^2}\:dx\right\}.\tag{2.2}
\label{2.2}
\end{equation*}
Furthermore, by applying the substitution $\frac{1-x}{1-ix}\mapsto x$ to the rightmost expression and using the identity $\int _0^1\frac{z\ln ^3\left(x\right)}{1-zx}\:dx=-6\operatorname{Li}_4\left(z\right)$ for $z\in \mathbb{C}\setminus \left(1,\infty \right)$, which is derived in more general terms in \cite[pp. 35–36]{book1} and \cite[pp. 70–71]{book2}, it follows that
$$\int _0^1\frac{\ln ^3\left(\frac{1-x}{1-ix}\right)}{1+x^2}\:dx=\frac{1}{2}\int _0^1\frac{\left(1-i\right)\ln ^3\left(x\right)}{1-\left(1+i\right)x}\:dx=-\frac{1}{2}i\int _0^1\frac{\left(1+i\right)\ln ^3\left(x\right)}{1-\left(1+i\right)x}\:dx=3i\operatorname{Li}_4\left(1+i\right).$$
Hence, substituting this into $\left(\ref{2.2}\right)$ and using the fact that $\int _0^1\frac{\arctan ^3\left(x\right)}{1+x^2}\:dx=\frac{45}{512}\zeta \left(4\right)$, we arrive at
$$\int _0^{\frac{\pi }{2}}x\ln ^2\left(1-\sin \left(x\right)\right)\:dx=\frac{15}{32}\zeta \left(4\right)+16\operatorname{\mathfrak{I}} \left\{i\operatorname{Li}_4\left(1+i\right)\right\}$$
$$=\frac{15}{32}\zeta \left(4\right)+16\operatorname{\mathfrak{R}} \left\{\operatorname{Li}_4\left(1+i\right)\right\}.$$
Consequently, applying the result in point $\left(\operatorname{i}\right)$ yields the corresponding closed-form expression.
\end{proof}
\end{lemma}

\begin{lemma} \label{lma4}Let $k$ be a positive integer. Then, the following equality holds:
$$\int _0^{\frac{\pi }{2}}x^2\cos ^{2k-1}\left(x\right)\:dx=\frac{1}{4}\frac{4^kH_k^{\left(2\right)}}{k\binom{2k}{k}}-\frac{4^kH_{2k}^{\left(2\right)}}{k\binom{2k}{k}}+\frac{3}{4}\zeta \left(2\right)\frac{4^k}{k\binom{2k}{k}}.$$
\begin{proof}
To begin, denote the integral of interest by $J_k$ and recall a variation of Wallis' integral, namely
$$I_n=\int _0^{\frac{\pi }{2}}\cos ^{2n-1}\left(x\right)\:dx=\frac{1}{2}\frac{4^n}{n\binom{2n}{n}}.$$
By applying integration by parts twice, we obtain
$$I_n=\frac{1}{2}\left(2n-1\right)\left(2n-2\right)\int _0^{\frac{\pi }{2}}x^2\cos ^{2n-3}\left(x\right)\:dx-\frac{1}{2}\left(2n-1\right)^2\int _0^{\frac{\pi }{2}}x^2\cos ^{2n-1}\left(x\right)\:dx,$$
or simply,
$$\frac{1}{2}\frac{4^n}{n\binom{2n}{n}}=\frac{1}{2}\left(2n-1\right)\left(2n-2\right)J_{n-1}-\frac{1}{2}\left(2n-1\right)^2J_n.$$
Now, if we multiply both sides by $\frac{2n\binom{2n}{n}}{4^n}\cdot \frac{1}{\left(2n-1\right)^2}$, we get
$$\frac{1}{\left(2n-1\right)^2}=\frac{2n-2}{2n-1}\frac{n\binom{2n}{n}}{4^n}J_{n-1}-\frac{n\binom{2n}{n}}{4^n}J_n,$$
and, by noting that $\binom{2n}{n}=\frac{2\left(2n-1\right)}{n}\binom{2n-2}{n-1}$ and summing both sides from $n=2$ to $k$, we arrive at
$$\sum _{n=2}^k\frac{1}{\left(2n-1\right)^2}=\sum _{n=2}^k\left(\frac{\left(n-1\right)\binom{2n-2}{n-1}}{4^{n-1}}J_{n-1}-\frac{n\binom{2n}{n}}{4^n}J_n\right).$$
Thus, by applying the identity $\sum _{n=2}^k\frac{1}{\left(2n-1\right)^2}=H_{2k}^{\left(2\right)}-\frac{1}{4}H_k^{\left(2\right)}-1$ and noting that the right-hand side forms a telescoping sum, it follows that
$$H_{2k}^{\left(2\right)}-\frac{1}{4}H_k^{\left(2\right)}-1=\frac{1}{2}J_1-\frac{k\binom{2k}{k}}{4^k}J_k.$$
Therefore,
$$H_{2k}^{\left(2\right)}-\frac{1}{4}H_k^{\left(2\right)}-1=\frac{1}{2}\int _0^{\frac{\pi }{2}}x^2\cos \left(x\right)\:dx-\frac{k\binom{2k}{k}}{4^k}\int _0^{\frac{\pi }{2}}x^2\cos ^{2k-1}\left(x\right)\:dx.$$
Hence, noting that $\int _0^{\frac{\pi }{2}}x^2\cos \left(x\right)\:dx=\frac{3}{2}\zeta \left(2\right)-2$ and isolating the desired integral, the identity follows.
\end{proof}
\end{lemma}

\section{The Main Series and their Proofs}
\begin{theorem}The following equalities hold:
\begin{alignat*}{3}
\left(\operatorname{i}\right)&\quad&&\sum _{k=1}^{\infty }\frac{4^kH_{2k}^{\left(2\right)}}{k^2\binom{2k}{k}}&&=-8G^2+11\zeta \left(4\right)+2\operatorname{Li}_4\left(\frac{1}{2}\right)+\ln ^2\left(2\right)\zeta \left(2\right)+\frac{1}{12}\ln ^4\left(2\right);\\
\left(\operatorname{ii}\right)&\quad&&\sum _{k=1}^{\infty }\frac{4^kH_{2k}^2}{k^2\binom{2k}{k}}&&=8G^2+\frac{103}{2}\zeta \left(4\right)-22\operatorname{Li}_4\left(\frac{1}{2}\right)+7\ln ^2\left(2\right)\zeta \left(2\right)-\frac{11}{12}\ln ^4\left(2\right).
\end{alignat*}
\end{theorem}

\begin{proof}
To evaluate the series in point $\left(\operatorname{i}\right)$, we begin by applying the identity in Lemma \ref{lma4}. If we multiply it by $\frac{1}{k}$, sum both sides from $k=1$ to $\infty$, use that $\sum _{k=1}^{\infty }\frac{x^{2k}}{k}=-\ln \left(1-x^2\right)$ for $\left|x\right|<1$, and substitute $x$ for $\cos \left(x\right)$, we obtain
$$-2\int _0^{\frac{\pi }{2}}x^2\sec \left(x\right)\ln \left(\sin \left(x\right)\right)\:dx=\frac{1}{4}\sum _{k=1}^{\infty }\frac{4^kH_k^{\left(2\right)}}{k^2\binom{2k}{k}}-\sum _{k=1}^{\infty }\frac{4^kH_{2k}^{\left(2\right)}}{k^2\binom{2k}{k}}+\frac{3}{4}\zeta \left(2\right)\sum _{k=1}^{\infty }\frac{4^k}{k^2\binom{2k}{k}}.$$
Moreover, by applying the results in point $\left(\operatorname{iii}\right)$ of Lemma \ref{lma1} and point $\left(\operatorname{iii}\right)$ of Lemma \ref{lma2}, we get
$$\frac{1}{4}\sum _{k=1}^{\infty }\frac{4^kH_k^{\left(2\right)}}{k^2\binom{2k}{k}}-\sum _{k=1}^{\infty }\frac{4^kH_{2k}^{\left(2\right)}}{k^2\binom{2k}{k}}=8G^2-\frac{45}{4}\zeta \left(4\right).$$
Thus, by isolating the target series and applying the result in point $\left(\operatorname{iv}\right)$ of Lemma \ref{lma1}, we obtain the desired closed form.
\\[1.5ex]
To carry out the evaluation of the series in point $\left(\operatorname{ii}\right)$, consider the result in point $\left(\operatorname{i}\right)$ of Lemma \ref{lma1}. If we let $k\mapsto 2k$, multiply both sides of the identity by $2\cdot \frac{4^k}{k\binom{2k}{k}}$, sum over $k=1$ to $\infty $, and apply the generating function in point $\left(\operatorname{ii}\right)$ of Lemma \ref{lma1}, it follows that
$$4\int _0^1\frac{\arcsin \left(x\right)\ln ^2\left(1-x\right)}{\sqrt{1-x^2}}\:dx=4\int _0^{\frac{\pi }{2}}x\ln ^2\left(1-\sin \left(x\right)\right)\:dx$$
$$=\sum _{k=1}^{\infty }\frac{4^kH_{2k}^2}{k^2\binom{2k}{k}}+\sum _{k=1}^{\infty }\frac{4^kH_{2k}^{\left(2\right)}}{k^2\binom{2k}{k}}.$$
Therefore, by isolating the corresponding series, applying the result of the integral in point $\left(\operatorname{ii}\right)$ of Lemma \ref{lma3}, and using the closed form of the series in point $\left(\operatorname{i}\right)$, we thereby obtain the closed-form expression.
\end{proof}

\end{document}